\documentclass[preprint,11pt]{elsarticle}

\newtheorem{theorem}{Theorem}[section]
\newtheorem{corollary}{Corollary}[section]
\newtheorem{lemma}{Lemma}[section]

\usepackage{amssymb}

\begin{document}

\begin{frontmatter}

\title{Exact rates of convergence in some martingale central limit theorems}
\author{Xiequan Fan}
 \cortext[cor1]{\noindent Corresponding author. \\
\mbox{\ \ \ \ }\textit{E-mail}: fanxiequan@hotmail.com (X. Fan). }
\address{Center for Applied Mathematics,
Tianjin University, 300072 Tianjin,  China;}
\address{Regularity Team, Inria, France}

\begin{abstract}
Renz \cite{R96}, Ouchti \cite{O05}, El Machkouri and Ouchti \cite{EO07} and Mourrat \cite{M13} have established some tight bounds on the rate of convergence in the central limit theorem for  martingales. In the present paper a
 modification of the methods,  developed by Bolthausen \cite{B82} and Grama and Haeusler \cite{GH00},  is  applied
for obtaining exact rates of convergence in the central limit theorem for martingales with differences
having conditional moments of order $2+\rho, \rho>0$.   Our results generalise and strengthen
the bounds mentioned above.
\end{abstract}

\begin{keyword}
Martingales,  Central limit theorem,  Berry-Esseen bounds
\vspace{0.3cm}
\MSC Primary 60G42; 60F05;  Secondary 60E15
\end{keyword}

\end{frontmatter}


\section{Introduction}
Assume that we are given a sequence of martingale differences $(\xi _i,\mathcal{F}_i)_{i=0,...,n}, $ defined on some
 probability space $(\Omega ,\mathcal{F},\mathbf{P})$,  where $\xi
_0=0 $ and $\{\emptyset, \Omega\}=\mathcal{F}_0\subseteq ...\subseteq \mathcal{F}_n\subseteq
\mathcal{F}$ are increasing $\sigma$-fields. Set
\begin{equation}\label{matingal}
X_{0}=0,\ \ \ \ \ \ \ X_k=\sum_{i=1}^k\xi _i,\quad k=1,...,n.
\end{equation}
Then $X=(X_k,\mathcal{F}_k)_{k=0,...,n}$ is a martingale.
Let $\left\langle X\right\rangle $ be its conditional variance:
\begin{equation}\label{quad}
\left\langle X\right\rangle _0=0,\ \ \ \ \ \ \ \ \left\langle X\right\rangle _k=\sum_{i=1}^k\mathbf{E}[\xi _i^2|\mathcal{F}
_{i-1}],\quad k=1,...,n.
\end{equation}
Define
$$D(X_n)= \sup_{ x \in \mathbf{R}}\Big|\mathbf{P}(X_n \leq x)-\Phi \left( x\right) \Big|,$$
where $\Phi \left( x\right)$ is the distribution function of the standard normal random variable.
Denote by $\stackrel{\   \mathbf{P} \ }{\longrightarrow}$ convergence in probability.
According to the basic results of martingale central limit theory (see  the monograph Hall and Heyde \cite{HH80}), the ``conditional Lindeberg condition"
\begin{eqnarray*}
\sum_{i=1 }^n \mathbf{E}[ \xi_i^2\mathbf{1}_{\{|\xi_i| \geq \varepsilon \}} | \mathcal{F}_{i-1}] \stackrel{\   \mathbf{P} \ }{\longrightarrow} 0, \ \ \ \ \ \textrm{as}\  n \rightarrow \infty \ \textrm{for each} \ \varepsilon>0,
\end{eqnarray*}
and the ``conditional normalizing condition"
\begin{eqnarray*}
 \left\langle X\right\rangle _n  \stackrel{\   \mathbf{P}\   }{\longrightarrow} 1, \ \ \ \ \ \textrm{as}\  n \rightarrow \infty,
\end{eqnarray*}
together implies that $$D(X_n)\longrightarrow 0, \ \ \ \  \ \textrm{as} \ \ n \rightarrow \infty.$$
In this paper we are interested in bounds of the speed of convergence in  central limit theorem, usually termed ``Berry-Esseen bounds''.

For general martingales, we first recall the following Berry-Esseen bound due to Heyde and Brown \cite{HB70}.
For $1< p \leq 2,$ Heyde and Brown
proved that
\begin{equation}\label{hb02}
D(X_n) \leq C_p\Big(\mathbf{E}\big[ \big| \langle X\rangle_n-1\big|^p \big] +\sum_{i=1}^n \mathbf{E}[ |\xi_i|^{2p} ] \Big)^{1/(2p+1) },
\end{equation}
 where $C_p$ depends only on $p$.
 The proof of Heyde and Brown is based on the martingale version of the Skorokhod embedding scheme.
 This method seems to be unsuited to obtain (\ref{hb02}) for $p>2.$
Using a method developed by Bolthausen \cite{B82}, Haeusler \cite{H88}  gave an extension of   (\ref{hb02}) to all $p>1$. See also Joos \cite{J93}.
Moreover, Haeusler also gave an example to show  that the bound (\ref{hb02}) is optimal under the stated condition, that
is there exists a sequence of martingale differences $(\xi_k,\mathcal{F}_k)_{k\geq0},$ such that for all
$n$ large enough,
\begin{eqnarray*}
D(X_n)\Big(\mathbf{E}\big[ \big| \langle X\rangle_n-1\big|^p \big] +\sum_{i=1}^n \mathbf{E}[ |\xi_i|^{2p} ] \Big)^{-1/(2p+1) } \ \geq \  c_p,
\end{eqnarray*}
 where $c_p$ is a positive constant and does not  depend  on $n$.
For  more interesting Berry-Esseen bounds for martingales, we refer to  Dedecker and Merlev\`{e}de \cite{DM11}, where
the authors consider the rates of convergence for linear statistics $X_n=\sum_{i \in \mathbb{Z}} c_{n, i} \xi_i $  based on stationary martingale differences.
Their rates are (most of the time) optimal in term of Wasserstein distances. As an application,
using the comparison between the uniform distance and the Kantorovith distance, it leads
to a Berry-Esseen bound of order  $n^{-1/4} \sqrt{ \log n }$ when the $\xi_i$'s have a moment of order
$3$ (see the rate (1.8) of \cite{DM11}). This Berry-Esseen bound provides the best rate of
convergence (up to the $\sqrt{ \log n }$ term) under the stated condition. Indeed,  Bolthausen \cite{B82} gave a
counter-example showing that the rate $n^{-1/4}$ in the Berry-Esseen bound cannot be improved when
 $\xi_i$'s have finite moments  of order $3$.

However, for martingales having bounded differences,  the bound  (\ref{hb02}) is not the best possible. In fact,
an earlier result of Bolthausen \cite{B82} sates that if $|\xi_i|\leq \epsilon$ and $\langle X\rangle_n=1$ a.s., then
\begin{equation}\label{ffsqs}
D(X_n) \leq C \epsilon^3 n\log  n,
\end{equation}
where $C$ is a constant. Moreover, Bolthausen \cite{B82} also showed that there exists
a sequence of martingale  differences satisfying $|\xi_i|\leq 2/\sqrt{n}$ and $\langle X\rangle_n=1$ a.s., such that for all $n$ large enough,
\begin{equation} \label{ghklm}
D(X_n) \sqrt{n} /\log n \ \geq \  c ,
\end{equation}
where $c$  is a  positive constant and does not  depend on $n$. This means the bound (\ref{ffsqs}) is optimal in the case that $\epsilon$ is of order $1/\sqrt{n}.$
Relaxing the condition  $\langle X\rangle_n=1$ a.s., Bolthausen \cite{B82} then proved that  if $|\xi_i|\leq \epsilon$ a.s.,
then
 \begin{equation}\label{sfhmd}
D(X_n) \leq C  \Big( \epsilon^3 n\log  n +  \min\{ ||\langle X\rangle_n-1||_1^{1/3}, ||\langle X\rangle_n-1||_{\infty}^{1/2}  \} \Big).
\end{equation}
It seems that the term $||\langle X\rangle_n-1||_1^{1/3}$ in the last bound should be replaced by $||\langle X\rangle_n-1||_1^{1/3} + \epsilon^{2/3}$;
see Mourrat \cite{M13}. (Indeed, in the proof of Bolthausen's corollary, we found a   term $\gamma^2$ is missing for the estimation of $\mathbf{E}[(\widehat{S}-S)^2]$; see \cite{B82} for details.)

If $\mathbf{E}[\xi _i^2|\mathcal{F}
_{i-1}]=1/n$  and $\mathbf{E} [|\xi_{i}|^{2+\rho} | \mathcal{F}_{i-1} ] \leq 1/n^{1+\rho/2}$ a.s.\ for some number $\rho \in (0, 1]$ and all $i=1,...,n, $
Renz \cite{R96} has obtained the following  Berry-Esseen bound:
\begin{eqnarray} \label{fsfsfgqs}
\ \ \ \ \ \ \ \  D(X_n) \, \leq \, C_{\rho} \, \varepsilon_n,
\end{eqnarray}
where the constant $C_\rho$ depends only on $\rho$   and
\begin{displaymath}
\varepsilon_n = \left\{ \begin{array}{ll}
  n^{-\rho/2}, & \textrm{\ \ \ if $\rho \in (0, 1),$}\\
 n^{-1/2}\log n , & \textrm{\ \ \ if $\rho=1.$}
\end{array} \right.
\end{displaymath}
Moreover, Renz  also showed that there exists
a sequence of martingale  differences satisfying his conditions, such that for all $n$ large enough,
\begin{equation}
D(X_n)\varepsilon_n^{-1}   \ \geq \  c ,
\end{equation}
where $c$ is a positive constant and does not depend on $n$. This means the bound (\ref{fsfsfgqs}) is exact.

With Bolthausen's method, El Machkouri and Ouchti   \cite{EO07} improved the term $ \epsilon^3 n\log  n $ in (\ref{sfhmd}) to $\epsilon \log  n,$ that is  if $\mathbf{E} [|\xi_{i}|^{3 } | \mathcal{F}_{i-1} ] \leq \epsilon \, \mathbf{E} [ \xi_{i} ^{2} | \mathcal{F}_{i-1} ]$ a.s.,
then
 \begin{equation}\label{sdf}
D(X_n) \leq C  \Big( \epsilon \,\log  n +  ||\langle X\rangle_n-1||_{\infty}^{1/2}    \Big).
\end{equation}
They also proved a result with term $||\langle X\rangle_n-1||_1^{1/3}.$ 

 Following Bolthausen \cite{B82} again, Mourrat \cite{M13} has obtained that if $|\xi_i|\leq \epsilon$ a.s., then for $p\geq1$,
\begin{eqnarray} \label{itmhklm}
\ \ \ \ \ \ \ \  D(X_n) \, \leq \, C_{p} \,   \Big(\epsilon^3 n\log  n+ \epsilon^{2p/(2p+1)} +  \mathbf{E}[| \left\langle X\right\rangle_n -1|^p] ^{1/(2p+1) } \Big),
\end{eqnarray}
where   $C_{p} $ is a constant depending only on $p.$  Notice that
Mourrat \cite{M13} has extended the term $ \min\{ ||\langle X\rangle_n-1||_1^{1/3}, ||\langle X\rangle_n-1||_{\infty}^{1/2}\}$ of Bolthausen \cite{B82} to the more general  term
$\mathbf{E}[| \left\langle X\right\rangle_n -1|^p] ^{1/(2p+1) } + \epsilon^{2p/(2p+1)} .$ Moreover, he also has justified the optimality
of the   term $\mathbf{E}[| \left\langle X\right\rangle_n -1|^p] ^{1/(2p+1) }.$

In this paper we give an improvement on the inequality of El Machkouri and Ouchti (\ref{sdf}) and Mourrat's inequality (\ref{itmhklm}).
 Our result also generalises the inequality of  Renz (\ref{fsfsfgqs}).   With the method of Grama and Haeusler \cite{GH00}, we prove that if there exist two positive numbers $\rho$ and $\epsilon$, such that
\begin{eqnarray}\label{gfvom1}
 \mathbf{E} [|\xi_{i}|^{2+\rho} | \mathcal{F}_{i-1} ] \leq \epsilon^{ \rho} \mathbf{E} [\xi_{i}^{2} | \mathcal{F}_{i-1} ]   \ \ \textrm{a.s.\ for all} \ i=1,...,n,
\end{eqnarray}
then
\begin{eqnarray}\label{ghka}
\ \ \ \ \ \ \ \  D(X_n) \, \leq \, C_{\rho} \,   \Big(  \gamma + ||\langle X\rangle_n-1||_{\infty}^{1/2}  \Big),
\end{eqnarray}
where  $C_{\rho} $ is a constant depending only on $\rho$ and
\begin{displaymath}
\gamma = \left\{ \begin{array}{ll}
\epsilon^\rho, & \textrm{\ \ \ if $\rho \in (0, 1),$}\\
\epsilon \, |\log \epsilon |, & \textrm{\ \ \ if $\rho\geq1.$}
\end{array} \right.
\end{displaymath}
  We also justify the optimality
of the term  $\gamma.$
Then with the method of Bolthausen \cite{B82}, we  obtain a  significant improvement of  Mourrat's inequality  (\ref{itmhklm})
by dropping the term $\epsilon^3 n\log  n $: If $|\xi_i|\leq \epsilon$ a.s., then for any  $p\geq 1$,
\begin{eqnarray} \label{dfsf}
\ \ \ \ \ \ \ \  D(X_n) \, \leq \, C_{p} \,   \Big(  \epsilon^{2p/(2p+1)}+ \mathbf{E}[| \left\langle X\right\rangle_n -1|^p] ^{1/(2p+1) } \Big),
\end{eqnarray}
where   $C_{p} $ is a constant depending only on $p.$

The paper is organized as follows. Our main results are stated and discussed in Section \ref{sec2}.
Proofs are deferred to Section \ref{sec3}.

Throughout the paper, $c$ and $c_\alpha$ probably supplied with some indices,
denote respectively a generic positive absolute constant and a generic positive constant depending only on $\alpha.$

\section{Main Results}\label{sec2}
In the sequel we shall use the following conditions:
\begin{description}
\item[(A1)]   There exist two positive numbers $\rho $ and $\epsilon \in (0, \frac12]$, such that
for all $1\leq i\leq n,$
\[
 \mathbf{E} [|\xi_{i}|^{2+\rho} | \mathcal{F}_{i-1} ] \leq \epsilon^{ \rho}  \mathbf{E}[\xi _i^2|\mathcal{F}
_{i-1}]\ \ \textrm{a.s.};
\]
\item[(A2)]  There exists a number  $ \delta\in [0, \frac12]$, such that
$ \left| \left\langle X\right\rangle _n-1\right| \leq  \delta^2$ a.s.
\end{description}

Let us comment on conditions (A1) and (A2).
\begin{enumerate}
\item Note that in the case of normalized sums of i.i.d.\ random variables, conditions (A1) and (A2) are satisfied with $\epsilon = \frac {1} {\sigma \sqrt n}$ and $\delta = 0$.

\item In the case of martingales, $\epsilon$ and $\delta$ usually depend on $n$ such that $\epsilon=\epsilon_n \rightarrow 0$ and $\delta=\delta_n \rightarrow 0$ as $n\rightarrow \infty$. It is also worth noting  that the bounded differences, that is  $|\xi_i|\leq\epsilon$ a.s.\ for all $i$, satisfy  condition (A1).

\item Assume that $(Y_i, \mathcal{F}_{i}  )_{i\geq1}$ is a sequence of martingale differences satisfying
 $$ \mathbf{E} [|Y_{i}|^{2+\rho} | \mathcal{F}_{i-1} ] \leq C^{\rho} \mathbf{E} [Y_{i}^{2} | \mathcal{F}_{i-1} ]$$
 for a positive absolute constant $C$ and  all $ i \geq 1.$
Let $S_n= \sum_{i=1}^nY_i$ and $s_n=\sqrt{\mathbf{E}[S_n^2]}.$
Then it is easy to verify that condition (A1) is satisfied with  $\xi_i =Y_i/s_n$ and  $\epsilon =C/s_n$. In particular, if  $(Y_i, \mathcal{F}_{i}  )_{i\geq 1}$ is a stationary sequence, then  $\epsilon =O(1/\sqrt{n})$ as $n\rightarrow\infty$.

\item Condition (A1) is satisfied for separately Lipschitz functions of independent random variables.
Let $f: {\mathcal X}^n \mapsto {\mathbf R}$ be separately Lipschitz, such that
\begin{equation} \label{codiMD}
|f(x_1, x_2, \ldots, x_n)-f(x'_1, x'_2, \ldots, x'_n)| \leq d(x_1,x'_1)+d(x_2,x'_2)+\cdots +  d(x_n, x'_n) \, .
\end{equation}
Let then
\begin{equation}\label{Sn}
X_n:=f(\eta_1, \ldots , \eta_n) -{\mathbf E}[f(\eta_1, \ldots , \eta_n)]\, ,
\end{equation}
where $\eta_1, \ldots , \eta_n$ is a sequence of independent random variables.
We also introduce the natural filtration of the chain, that is  ${\mathcal F}_0=\{\emptyset, \Omega \}$
and for $k \in {\mathbf{N}}$,
${\mathcal F}_k= \sigma(\eta_1, \eta_2,  \ldots, \eta_k)$.
Define then
\begin{equation}\label{gk}
g_k(\eta_1, \ldots , \eta_k)= {\mathbf E}[f(\eta_1, \ldots, \eta_n)|{\mathcal F}_k]\, ,
\end{equation}
and
\begin{equation}\label{dk}
\xi_k=g_k(\eta_1, \ldots, \eta_k)-g_{k-1}(\eta_1, \ldots, \eta_{k-1})\, .
\end{equation}
For $k \in [1, n-1]$, let
$$
X_k:=\xi_1+\xi_2+\cdots + \xi_k \, ,
$$
and note that, by definition of the $\xi_k$'s, the functional $X_n$ introduced in (\ref{Sn})
 satisfies
$$
X_n =\xi_1+\xi_2+\cdots + \xi_n \, .
$$
Hence $X_k$ is a martingale adapted to the filtration ${\mathcal F}_k$.  It is easy to verify that
$\xi_1,...,\xi_n$ satisfy condition (A1). Indeed, for all $1\leq i\leq n,$
\begin{eqnarray}
 \mathbf{E} [|\xi_{i}|^{2+\rho} | \mathcal{F}_{i-1} ]&=&  \mathbf{E}[  |{\mathbf E}[f(\eta_1, \ldots, \eta_n)|{\mathcal F}_i] -{\mathbf E}[f(\eta_1, \ldots, \eta_n)|{\mathcal F}_{i-1}]|^{\rho}\xi_{i}^2 \, | \mathcal{F}_{i-1}]  \nonumber \\
 &=&  \mathbf{E}[  |{\mathbf E}[f(\eta_1, \ldots, \eta_n)|{\mathcal F}_i] -{\mathbf E}[f(\eta_1, \ldots, \eta_i', \ldots, \eta_n)|{\mathcal F}_{i}]|^{\rho}\xi_{i}^2 \, | \mathcal{F}_{i-1}]  \nonumber  \\
&\leq&  \mathbf{E}[  ({\mathbf E}[d(\eta_i, \eta_i')|{\mathcal F}_i])^{\rho}\xi_{i}^2 \, | \mathcal{F}_{i-1}] \nonumber  \\
&=& \mathbf{E}[({\mathbf E}[d(\eta_i, \eta_i')|\eta_i])^{\rho}]\  \mathbf{E}[  \xi_{i}^2 \, | \mathcal{F}_{i-1}],  \label{condition25}
\end{eqnarray}
where $(\eta_1', \ldots , \eta_n')$ is an independent copy of $(\eta_1, \ldots , \eta_n).$
Hence, condition (A1) is satisfied with $\epsilon =\max_{1\leq i \leq n} \mathbf{E}[({\mathbf E}[d(\eta_i, \eta_i')| \eta_i])^{\rho}]^{1/\rho}.$
In particular, by Jensen's inequality, it holds $\epsilon  \leq \max_{1\leq i \leq n} {\mathbf E}[d(\eta_i, \eta_i')]$
for $\rho \in (0, 1].$
\end{enumerate}

Our  first result is the following   Berry-Esseen bounds  for martingales.
\begin{theorem}\label{th0}
Assume conditions (A1) and (A2).
\begin{itemize}
  \item If $\rho \in (0, 1)$, then
\begin{equation}\label{t0ie1}
D(X_n)\leq c_{\rho} \Big( \epsilon^\rho   + \delta \Big) .
\end{equation}

  \item If $\rho \in [1, \infty)$, then
\begin{equation} \label{t1ie2}
D(X_n)\leq c \, \Big( \epsilon  |\log \epsilon|  + \delta \Big).
\end{equation}
\end{itemize}
\end{theorem}

We   justify the optimality
of the term  $\epsilon^\rho$ of (\ref{t0ie1}). Let $n= \lfloor\epsilon^{-2}\rfloor$  be the integer part of $\epsilon^{-2}$ and $\rho \in (0, 1)$.   Renz's
 inequality   (\ref{ghklm}) shows that there exists a sequence of martingale  differences satisfying condition (A1)
 and $\langle X\rangle_n=1$ a.s., such that for all $\epsilon$ small enough,
\begin{equation}
\epsilon^{-\rho} D(X_n) \geq n^{\rho/2} D(X_n) \geq    c ,
\end{equation}
where the constant $c>0$   does not depend on $\epsilon$.

Notice that, for  bounded martingale differences,  condition (A1) holds with   $\rho=1$.
By Bolthausen's inequality (\ref{ghklm}) with $n=[\epsilon^{-2}]$, there exists
a sequence of martingale  differences satisfying $|\xi_i| \leq 3\epsilon$ and $\langle X\rangle_n=1$ a.s., such that for all $\epsilon$ small enough,
\begin{equation}
(3\epsilon \, |\log 3\epsilon | )^{-1} D(X_n) \geq \frac14 D(X_n)\sqrt{n} /\log n \geq    c ,
\end{equation}
where the constant $c>0$   does not depend on $\epsilon$. Thus the term  $\epsilon \, |\log \epsilon |$ of (\ref{t1ie2}) is exact even for bounded martingale differences.

Under the conditions (A1) and (A2), the order of the term $\epsilon \left| \log  \epsilon\right|$ in (\ref{t1ie2}) is less than the order of the term $\epsilon^3 n \log  n  $ in Bolthausen's inequality (\ref{sfhmd}). Indeed, by condition (A2), we have $3/4 \leq \langle X \rangle_n \leq n \epsilon ^2$ a.s.\ (see Lemma \ref{lemma2s}) and then $\epsilon \geq \sqrt{3/(4n)}$. For $\epsilon \leq 1/2$, it is easy to see that $\epsilon^3 n   \log n  \geq 3\,\epsilon|\log\epsilon | /4$.
Moreover, $\epsilon^3 n \log  n $ may converge to infinity while $\epsilon \left| \log  \epsilon\right|$  converges to $0$ as $\epsilon  \to 0$ and $n \to \infty.$  For instance, if $\epsilon$ is of the order $n^{-1/3}$  as $n \to \infty$, then it is obvious that $\epsilon \left| \log  \epsilon \right|=O(n^{-1/3} \log n  )$ while $\epsilon^3 n  \log n  \geq \log n.$ Thus the term  $\epsilon \left| \log  \epsilon\right|$
 is much smaller than $\epsilon^3 n \log  n  $.
Similarly,  the order of $ \epsilon\,|\log \epsilon |$ is also better than the order of $\epsilon\,\log n$ in (\ref{sdf}) of El Machkouri and Ouchti   \cite{EO07}.

For martingales with bounded differences, inequality (\ref{t1ie2}) has been established earlier in Grama
\cite{G87a,G87b}. Under the conditional Bernstein condition, that is
\[
|\mathbf{E}[\xi_{i}^{k}  | \mathcal{F}_{i-1}]| \leq \frac12 k!\epsilon^{k-2} \mathbf{E}[\xi_{i}^2 | \mathcal{F}_{i-1}]\ \   \textrm{a.s.\ for}\ k\geq 3\ \ \textrm{and}\ \ 1\leq i\leq n,
\]
 instead of condition (A1), Fan, Grama and Liu \cite{F13} have
obtained the Berry-Esseen bound (\ref{t1ie2}).  Note that the conditional Bernstein condition
implies that  $\xi_i$ has conditional exponential moment.
Now we only  assume that  $\xi_i$ has conditional moment of order $3.$

Using Theorem  \ref{th1}, we have the following Berry-Esseen bounds similar to the results of Ouchti \cite{O05}. Following the notations of Ouchti \cite{O05},
let $v(n)$ denote either
$$ \sup\{k:\ \langle X\rangle_k \leq1\}   \ \ \ \
\textrm{or } \ \ \
  \inf\{k:\ \langle X\rangle_k \geq  1\}.$$
\begin{corollary}\label{co20}
Assume conditions (A1) and   $\langle X\rangle_n \geq  1$  a.s.
\begin{itemize}
  \item If $\rho \in (0, 1)$, then
\begin{equation}\label{tjhk1}
D(X_{v(n)} ) \leq c_{\rho} \,  \epsilon^\rho .
\end{equation}

  \item If $\rho \in [1, \infty)$, then
\begin{equation} \label{tjhk2}
D(X_{v(n)} ) \leq c \,  \epsilon  |\log \epsilon| .
\end{equation}
\end{itemize}
\end{corollary}

Inequality (\ref{tjhk2}) significantly improves an earlier result of Ouchti \cite{O05} under the following condition $$\mathbf{E} [|\xi_{i}|^{3} | \mathcal{F}_{i-1} ] \leq n^{-1/2} \mathbf{E}[\xi _i^2|\mathcal{F}
_{i-1}] \ \ \textrm{a.s. for all}\ i\geq 1. $$ Ouchti has obtained a convergence rate in  central limit theorem  of order $n^{-1/4},$
while (\ref{tjhk2}) gives a convergence rate
of order $n^{-1/2}\log n.$

Relaxing  condition  (A2), we have the following estimation.
\begin{theorem}\label{th1}
Assume condition  (A1). Let  $p\geq 1$.
\begin{itemize}
  \item If $\rho \in (0, 1)$, then
\begin{equation}\label{sffs}
D(X_n)\leq c_{p,\rho} \, \bigg( \epsilon^\rho+\Big(\mathbf{E}\big[ \big| \langle X\rangle_n-1\big|^p \big] +\mathbf{E}[ \max_{1 \leq i \leq n} |\xi_i|^{2p} ]\Big)^{1/(2p+1) } \bigg) .
\end{equation}

  \item If $\rho \in [1, \infty)$, then
\begin{equation}\label{sff026ss}
D(X_n)\leq c_{p} \,   \bigg( \epsilon |\log \epsilon | +\Big(\mathbf{E}\big[ \big| \langle X\rangle_n-1\big|^p \big] +\mathbf{E}[ \max_{1 \leq i \leq n} |\xi_i|^{2p} ]\Big)^{1/(2p+1) } \bigg) .
\end{equation}
\end{itemize}
\end{theorem}

Notice that $ \mathbf{E}[ \max_{1 \leq i \leq n} |\xi_i|^{2p} ]  \leq   \sum_{i=1}^n \mathbf{E}[  |\xi_i|^{2p} ] .$
Therefore, our bounds are  usually smaller than the bound of Haeusler   (\ref{hb02}). For instance, if $|\xi_i|\leq \epsilon$ a.s.\
for all $i=1,...,n,$ then
$ \mathbf{E}[ \max_{1 \leq i \leq n} |\xi_i|^{2p} ] \leq \epsilon^{2p}$,   while $\sum_{i=1}^n \mathbf{E}[  |\xi_i|^{2p} ] \leq n\epsilon^{2p}.$

%
%

For martingales having bounded differences, Theorem \ref{th1}
implies the following corollary.
\begin{corollary}\label{co1}
Assume $| \xi_i|\leq \epsilon$ a.s.\ for all $i \in [0, n]$.
 Then  for any  $p\geq 1$,
\begin{eqnarray}\label{sfnjm2}
D(X_n) \leq c_p \Big( \epsilon^{2p} + \mathbf{E}\big[\big| \langle X\rangle_n-1\big|^p \big]\Big)^{1/(2p+1)}.
\end{eqnarray}
\end{corollary}

Clearly,  the term $\epsilon^3 n\log  n$ appearing in Mourrat's inequality  (\ref{itmhklm}) does not appear any more in (\ref{sfnjm2}).
When $\epsilon \rightarrow 0$ and $\epsilon \geq \sqrt[3/7]{1/(n \log n)},$  it holds  $\epsilon^{2p/(2p+1)} \leq \epsilon^3 n\log  n$ for any $p\geq 1$.
Moreover, when $\epsilon \rightarrow 0$ and $\epsilon \geq 1/ \sqrt[3]{n},$ we have $\epsilon^3 n\log  n \rightarrow \infty$ and  $\epsilon^{2p/(2p+1)}\rightarrow0.$
Thus our bound (\ref{sfnjm2}) is
significantly smaller than the bound of  Mourrat  (\ref{itmhklm}).

\section{Proofs of Theorems}\label{sec3}
In the sequel, for simplicity, the equalities and inequalities involving random variables will be understood in the a.s.\ sense without mentioning this.

In the proofs of theorems, we will make use of the following two lemmas.
The first lemma shows  that we may
assume  $\rho \in (0, 1]$ in condition (A1).
\begin{lemma}\label{lema1} If there  exists an $s> 2,$ such that
\begin{equation}\label{ineq31}
\mathbf{E} [|\xi_{i}|^{s} | \mathcal{F}_{i-1} ] \leq \epsilon^{s-2} \, \mathbf{E} [\xi_{i}^{2} | \mathcal{F}_{i-1} ],
\end{equation}
then, for any $ t \in [2,   s)$,
\begin{equation}
\mathbf{E} [|\xi_{i}|^{t} | \mathcal{F}_{i-1} ] \leq \epsilon^{t-2} \, \mathbf{E} [\xi_{i}^{2} | \mathcal{F}_{i-1} ].
\end{equation}
\end{lemma}
\emph{Proof.}
Let $l, p, q$ be defined by the following equations
 $$lp=2, \ \ \ \ (t-l)q=s, \ \ \  p^{-1}+ q^{-1}=1, \ \ \ l>0\ \textrm{and}\  p, q\geq1.$$
Solving the last equations, we get
 $$l= \frac{2(s-t)}{s-2},\ \ \ \ \  p=\frac{s-2}{s-t}, \ \ \ \ \  q=\frac{s-2}{t-2}.$$
By  H\"{o}lder's inequality and (\ref{ineq31}), it is easy to see that
\begin{eqnarray*}
\mathbf{E} [|\xi_{i}|^{t} | \mathcal{F}_{i-1} ] &=& \mathbf{E} [|\xi_{i}|^{l} |\xi_{i}|^{t-l} | \mathcal{F}_{i-1} ] \\
&\leq &( \mathbf{E} [|\xi_{i}|^{lp}  | \mathcal{F}_{i-1} ])^{1/p} ( \mathbf{E} [  |\xi_{i}|^{(t-l)q} | \mathcal{F}_{i-1} ])^{1/q}\\
&\leq & ( \mathbf{E} [\xi_{i}^{2}  | \mathcal{F}_{i-1} ])^{1/p} ( \mathbf{E} [  |\xi_{i}|^{s} | \mathcal{F}_{i-1} ])^{1/q} \\
&\leq & ( \mathbf{E} [\xi_{i}^{2}  | \mathcal{F}_{i-1} ])^{1/p} (\epsilon^{s-2} \mathbf{E} [ \xi_{i}^2 | \mathcal{F}_{i-1} ])^{1/q} \\
&\leq & \epsilon^{(s-2)/q} \mathbf{E} [\xi_{i}^{2}  | \mathcal{F}_{i-1} ]\\
&=& \epsilon^{t-2} \ \mathbf{E} [\xi_{i}^{2} | \mathcal{F}_{i-1} ].
\end{eqnarray*}
This completes the proof of lemma.  \qed

The following lemma shows that
under condition (A1), $\xi_{i}$ has a bounded conditional variance.
\begin{lemma}\label{lemma2s}
If there exists an $s> 2,$  such that
\begin{equation}
\mathbf{E} [|\xi_{i}|^{s} | \mathcal{F}_{i-1} ] \leq \epsilon^{s-2} \, \mathbf{E} [\xi_{i}^{2} | \mathcal{F}_{i-1} ],
\end{equation}
then
\begin{equation}\label{ineq28}
\mathbf{E} [\xi_{i}^{2} | \mathcal{F}_{i-1} ] \leq \epsilon^2.
\end{equation}
In particular, condition (A1) implies    (\ref{ineq28}).
\end{lemma}
\emph{Proof.}
By Jensen's inequality, it is easy to see that
\begin{eqnarray*}
(\mathbf{E} [\xi_{i}^{2} | \mathcal{F}_{i-1} ])^{s/2} &\leq& \mathbf{E} [|\xi_{i}|^{s} | \mathcal{F}_{i-1} ] \\
&\leq & \epsilon^{s-2} \ \mathbf{E} [\xi_{i}^{2} | \mathcal{F}_{i-1} ].
\end{eqnarray*}
Thus $$(\mathbf{E} [\xi_{i}^{2} | \mathcal{F}_{i-1} ])^{s/2 -1} \leq   \epsilon^{s-2},$$
 which implies (\ref{ineq28}).  \qed

\subsection{Proof of Theorem \ref{th0}}
Theorem  \ref{th0} is a refinement of Lemma 3.3 of Grama and Haeusler \cite{GH00} where it is assumed that $\xi_{i}$'s are bounded,  which is a  particular case of condition  (A1). See also Lemma 3.1 of Fan, Grama and Liu \cite{F13}.   Compared to the  proofs of Grama and Haeusler \cite{GH00} and Fan, Grama and Liu \cite{F13},   the main challenge of our proof comes from  the control of $I_1$ defined in (\ref{D-1}).

By Lemma \ref{lema1},  we only need to consider the case of $\rho \in (0, 1].$
Set $T=1+\delta ^2$, and introduce a modification of the
conditional variance $\left\langle X\right\rangle $ as follows:
\begin{equation}
V_k=\left\langle X\right\rangle _k\mathbf{1}_{\{k<n\}}+T\mathbf{1}_{\{k=n\}}.
\label{RA-4}
\end{equation}
It is obvious that $V_0=0,$ $V_n=T$, and that $(V_k,\mathcal{F}_k)_{k=0,...,n}$ is a
predictable process. For simplicity of notations, denote
\begin{displaymath}
\gamma  = \left\{ \begin{array}{ll}
\epsilon  +\delta, & \textrm{\ \ \ if $\rho \in (0, 1),$}\\
\epsilon\left| \log   \epsilon  \right|+ \delta , & \textrm{\ \ \ if $\rho=1.$}
\end{array} \right.
\end{displaymath}
Let $c_{*}$ be a  constant depending on $\rho$, whose   value will be
chosen later. Define the following non-increasing discrete time predictable process
$$A_k=c_{*}^2\gamma ^2+T-V_k,\ \ \ \ \  k=1,...,n.$$
In particular, we have $A_0=c_{*}^2\gamma ^2+T$ and $A_n=c_{*}^2\gamma ^2$.
Moreover, for  $u, x\in \mathbf{R}$,  and $y > 0,$ set, for brevity,
\begin{equation}
\Phi _u(x,y)=\Phi \Big( \frac{u-x}{\sqrt{y}}  \Big) .  \label{RA-7}
\end{equation}

Let $\mathcal{N} =\mathcal{N}(0,1)$ be a standard normal random variable, which is independent of $X_n$.
Using a  smoothing procedure, by Lemma \ref{LEMMA-APX-1},
we get
\begin{eqnarray}
 \sup_u\Big| \mathbf{P} (X_n\leq u)  -  \Phi (u)\Big| &\leq& c_1\sup_u\Big|
 \mathbf{P} (X_n + c_{*} \gamma\mathcal{N}\leq u)  -  \Phi (u)\Big|
+c_2\gamma \nonumber\\
&=& c_1\sup_u\Big|
\mathbf{E} [\Phi _u(X_n,A_n)]- \Phi (u)\Big|
+c_2\gamma \nonumber\\
&\leq& c_1\sup_u\Big| \mathbf{E} [\Phi _u(X_n,A_n)]-\mathbf{E} [\Phi _u(X_0,A_0)]\Big| \nonumber\\
& & +\, c_1\sup_u\Big| \mathbf{E} [\Phi _u(X_0,A_0)]-\Phi (u)\Big| +c_2\gamma   \nonumber\\
&=& c_1\sup_u\Big| \mathbf{E} [\Phi _u(X_n,A_n)]-\mathbf{E} [\Phi _u(X_0,A_0)]\Big| \nonumber\\
& & +\, c_1\sup_u\Big| \Phi \Big(\frac{u}{\sqrt{c_{*}^2\gamma^2+T}}  \Big)-\Phi (u)\Big| +c_2\gamma .   \label{llaa01}
\end{eqnarray}
Since $T=1+\delta^2,$ it is easy to see that
\begin{eqnarray}
 \Big| \Phi \Big(\frac{u}{\sqrt{c_{*}^2\gamma^2+T}}  \Big)-\Phi (u)\Big| \leq  c_3 \Big|\frac{1}{\sqrt{c_{*}^2\gamma^2+T}}  -1 \Big|
\leq c_4\gamma .
\end{eqnarray}
Returning to (\ref{llaa01}), we obtain
\begin{eqnarray}
 \sup_u\Big| \mathbf{P} (X_n\leq u)  -  \Phi (u)\Big|  \leq
  c_1\sup_u\Big| \mathbf{E} [\Phi _u(X_n,A_n)]-\mathbf{E} [\Phi _u(X_0,A_0)]\Big|   +c_5\gamma .   \label{llaa022}
\end{eqnarray}
By a simple telescoping, we deduce that
\[
\mathbf{E} [\Phi _u(X_n,A_n)]-\mathbf{E} [\Phi _u(X_0,A_0)]=\mathbf{E}
\Big[ \sum_{k=1}^n \Big(\Phi _u(X_k,A_k)-\Phi _u(X_{k-1},A_{k-1}) \Big)\Big] .
\]
Using the fact
\[
\frac{\partial ^2}{\partial x^2}\Phi _u(x,y)=2\frac \partial {\partial
y}\Phi _u(x,y),
\]
we obtain
\begin{equation}
\mathbf{E}[ \Phi _u(X_n,A_n)]-\mathbf{E} [\Phi _u(X_0,A_0)]=I_1+I_2-I_3,
\label{lla}
\end{equation}
where
\begin{eqnarray}
I_1&=&\mathbf{E} \bigg[ \sum_{k=1}^n \bigg(\frac{}{} \Phi _u(X_k,A_k)-\Phi
_u(X_{k-1},A_k)  \nonumber \\
&& \ \ \ \  \ -\, \frac \partial {\partial x}\Phi
_u(X_{k-1},A_k)\xi_k-\frac 12\frac{\partial ^2}{\partial x^2}\Phi
_u(X_{k-1},A_k)\xi _k^2 \bigg) \bigg] ,  \label{D-1}\\
I_2 &=& \frac 12\mathbf{E}   \bigg[  \sum_{k=1}^n\frac{\partial ^2}{\partial x^2}\Phi
_u(X_{k-1},A_k)\Big(\Delta \left\langle X\right\rangle _k-\Delta V_k \Big) \bigg] ,\quad \quad \
\label{D-2} \\
I_3&=&\mathbf{E}   \bigg[\sum_{k=1}^n\left( \Phi _u(X_{k-1},A_{k-1})-\Phi
_u(X_{k-1},A_k)-\frac \partial {\partial y}\Phi _u(X_{k-1},A_k)\Delta
V_k\right)  \bigg],  \label{D-3}
\end{eqnarray}
where $\Delta \left\langle X\right\rangle _k= \left\langle X\right\rangle _k-\left\langle X\right\rangle _{k-1}$.

Next, we  give the estimates of $I_1,$ $I_2$ and $I_3.$ To this end, we introduce the following notations.
Denote by $\varphi$   the density function of the standard normal random variable.
Moreover, $\vartheta_i$'s stand for some values or random variables satisfying $0\leq \vartheta_i  \leq 1$,
which may represent different values at different places.

\emph{\textbf{a)} Control of} $I_1.$
To shorten notations, set
$T_{k-1}= (u-X_{k-1})/\sqrt{A_k}.$ It is obvious that
\begin{eqnarray}
&&  R_k=: \Phi _u(X_k,A_k)-\Phi
_u(X_{k-1},A_k)    -\, \frac \partial {\partial x}\Phi
_u(X_{k-1},A_k)\xi_k-\frac 12\frac{\partial ^2}{\partial x^2}\Phi
_u(X_{k-1},A_k)\xi _k^2 \nonumber  \\
&& \ \ \ \ \,    =\   \Phi \Big(T_{k-1} + \frac{\xi_k}{\sqrt{A_k}}\Big)-\Phi(T_{k-1})
    -\, \Phi'(T_{k-1})\frac{\xi_k}{\sqrt{A_k}}   -\frac 12\Phi''(T_{k-1})\Big(\frac{\xi_k}{\sqrt{A_k}}\Big)^2 . \nonumber
\end{eqnarray}
We distinguish two cases as follows.

\emph{Case 1}: $|\xi_{k}/  \sqrt{A_{k}}| \leq 1+   |T_{k-1}| /2$.
By a three-term Taylor expansion,   it is easy to see that  if  $|\xi_k/\sqrt{A_k}|\leq 1,$ then
\begin{eqnarray*}
\Big|R_k \Big|  & =&   \Big|\,\frac16 \Phi'''\Big(T_{k-1} + \vartheta  \frac{\xi_k}{\sqrt{A_k}}\Big)  \Big|\frac{\xi_k}{\sqrt{A_k}}\Big|^3 \, \Big| \nonumber \\
 &\leq&  \Big|  \Phi'''\Big(T_{k-1} + \vartheta  \frac{\xi_k}{\sqrt{A_k}}\Big)\Big|  \Big|\frac{\xi_k}{\sqrt{A_k}}\Big|^{2+\rho}.
\end{eqnarray*}
It is also easy to see that if  $|\xi_k/\sqrt{A_k}|> 1,$ then
\begin{eqnarray*}
\Big|R_k \Big|  &\leq & \frac12 \Big( \Big|  \Phi''\Big(T_{k-1} + \vartheta  \frac{\xi_k}{\sqrt{A_k}}\Big)\Big|  + \Big| \Phi''(T_{k-1}) \Big|  \Big) \, \Big(\frac{\xi_k}{\sqrt{A_k}}\Big)^{2} \nonumber \\
&\leq&  \Big|  \Phi''\Big(T_{k-1} + \vartheta'  \frac{\xi_k}{\sqrt{A_k}}\Big)\Big|  \, \Big(\frac{\xi_k}{\sqrt{A_k}}\Big)^{2}\nonumber \\
&\leq&   \Big|  \Phi''\Big(T_{k-1} + \vartheta'  \frac{\xi_k}{\sqrt{A_k}}\Big)\Big|  \Big|\frac{\xi_k}{\sqrt{A_k}}\Big|^{2+\rho} ,
\end{eqnarray*}
where
\begin{displaymath}
 \vartheta' = \left\{ \begin{array}{ll}
 \vartheta, & \textrm{\ \ \ if $\Big|  \Phi''\Big(T_{k-1} + \vartheta  \frac{\xi_k}{\sqrt{A_k}}\Big)\Big| \geq \Big| \Phi''(T_{k-1}) \Big|,$}\\
 \\
0 , & \textrm{\ \ \ if $ \Big|  \Phi''\Big(T_{k-1} + \vartheta  \frac{\xi_k}{\sqrt{A_k}}\Big)\Big|  < \Big| \Phi''(T_{k-1}) \Big|.$}
\end{array} \right.
\end{displaymath}
By the inequality $ \max\{ |\Phi''(t)|, \,  |\Phi'''(t)| \} \leq\varphi( t)(1+t^2)$, it follows that
\begin{eqnarray*}
   \Big|R_k\mathbf{1}_{\left\{|\xi_{k}/ \sqrt{A_{k}}| \leq  1+ |T_{k-1}| /2  \right\}} \Big| &\leq & \varphi\Big(T_{k-1}+\vartheta_1\frac{\xi_{k}}{ \sqrt{A_{k}}} \Big)\bigg(1+\Big(T_{k-1}+\vartheta_1\frac{\xi_{k}}{ \sqrt{A_{k}}} \Big)^2 \bigg) \\
  &\leq &  g_{1}(T_{k-1}) ,
\end{eqnarray*}
where $$g_{1}(z)=\sup_{|t-z|  \leq   1+ |z| /2 }\varphi(t )(1+ t^2).$$
It is easy to see that $g_{1}(z)$ is a  non-increasing in $z\geq 0$, and that $g_{1}(z)$ satisfies
\begin{eqnarray}\label{fklm59}
 \Big|R_k\mathbf{1}_{\left\{|\xi_{k}/ \sqrt{A_{k}}| \leq  1+ |T_{k-1}| /2  \right\}} \Big| \leq  g_{1}(T_{k-1}) \Big|\frac{\xi_{k}}{  \sqrt{A_{k}}} \Big|^{2+\rho} \mathbf{1}_{\left\{|\xi_{k}/ \sqrt{A_{k}}| \leq  1+ |T_{k-1}| /2  \right\}} .
\end{eqnarray}

\emph{Case 2}: $|\xi_{k}/ \sqrt{A_{k}}| > 1+ |T_{k-1}| /2$.
It is easy to see that for $|\Delta x|>1+   |x|/2,$
\begin{eqnarray*}
&&\Big|\Phi(x+\Delta x)-\Phi(x) - \Phi'(x)  \Delta x - \frac12 \Phi''(x) (\Delta x)^2 \Big|   \nonumber \\
& &=  \Big(  \Big| \frac{\Phi(x+\Delta x)-\Phi(x)}{|\Delta x|^{2+\rho}}  \Big|       +  | \Phi'(x)|   +  |   \Phi''(x)  |  \Big) \, |\Delta x|^{2+\rho} \nonumber \\
&&\leq   \Big(4 \Big| \frac{\Phi(x+\Delta x)-\Phi(x)}{(2+ |x|)^{2}}  \Big| +  | \Phi'(x)| + |\Phi''(x)|  \Big)  \, |\Delta x|^{2+\rho} \nonumber \\
&&\leq   \Big(  \frac{c_1}{(2+ |x|)^{2}}+  | \Phi'(x)| + |\Phi''(x)|  \Big)  \, |\Delta x|^{2+\rho}\nonumber \\
&&\leq   \frac{c_2}{(2+ |x|)^{2}}   \, |\Delta x|^{2+\rho}.
\end{eqnarray*}
Therefore,
\begin{eqnarray} \label{fklm60}
 \Big|R_k\mathbf{1}_{\left\{|\xi_{k}/ \sqrt{A_{k}}| > 1+ |T_{k-1}| /2  \right\}} \Big| \leq   g_{2}(T_{k-1}) \Big|\frac{\xi_{k}}{  \sqrt{A_{k}}} \Big|^{2+\rho} \mathbf{1}_{\left\{|\xi_{k}/ \sqrt{A_{k}}| > 1+ |T_{k-1}| /2  \right\}}  ,
\end{eqnarray}
where
$$g_2(z)=\frac{c_2}{(2+ |z|)^{2}} .$$

Set $$G(z)=  g_{1}(z)+g_{2}(z) .$$
Combining (\ref{fklm59}) and  (\ref{fklm60}) together, we obtain
\begin{eqnarray}
 \Big|R_k  \Big| \leq  G(T_{k-1}) \Big|\frac{\xi_{k}}{  \sqrt{A_{k}}} \Big|^{2+\rho}   ,
\end{eqnarray}
and thus
\begin{eqnarray}\label{ineq46}
 \Big|I_1 \Big| = \Big|\mathbf{E}  \Big[\sum_{k=1}^n  R_k   \Big] \Big| \leq \mathbf{E}  \Big[\sum_{k=1}^n G(T_{k-1}) \Big|\frac{\xi_{k}}{  \sqrt{A_{k}}} \Big|^{2+\rho}  \Big] .
\end{eqnarray}
Now we consider the conditional expectation of $|\xi_{k}|^{2+\rho}$. Using condition (A1), we have
\[
\mathbf{E}[|\xi_{k}|^{2+\rho} | \mathcal{F}_{k-1}] \leq  \epsilon^{ \rho} \,\Delta \langle X \rangle_{k},
\]
where $\Delta \langle X \rangle_{k}=\langle X \rangle_{k}-\langle X \rangle_{k-1}$.
It is obvious that $$\Delta \left\langle X\right\rangle _k=\Delta V_k=V_{k}-V_{k-1},\  \ 1\leq k <n, \ \
\Delta \left\langle X\right\rangle_n\leq \Delta V_n,$$ and that
\begin{equation} \label{ee}
\mathbf{E}[|\xi_{k}|^{2+\rho} | \mathcal{F}_{k-1}] \leq \epsilon^{ \rho}\, \Delta V_k.
\end{equation}
Combining (\ref{ineq46}) and  (\ref{ee}) together, we obtain
\begin{eqnarray}\label{fgfgfdgk}
\Big|I_1 \Big| \leq J_1:= \epsilon ^{\rho}\,\Big[\sum_{k=1}^n \frac{1}{ A_{k}^{1+\rho/2}} G(T_{k-1})\,  \Delta V_k   \Big] .
\end{eqnarray}

To estimate $J_1,$ we introduce the time change $\tau _t$ as follows:
for any real $t\in [0,T]$,
\begin{equation}
\tau _t=\min \{k\leq n: \, V_{k}>t\},\quad \textrm{where}\  \min \emptyset =n.
\end{equation}
Clearly, for any $t\in [0,T],$ the stopping time $\tau _t$ is
predictable. Denote by $(\sigma _k)_{k=1,...,n+1}$  the increasing sequence of
moments when the increasing and stepwise function $\tau _t,$ $t\in [0,T]$, has
jumps. It is obvious that $\Delta V_k =\int_{[\sigma _k,\sigma _{k+1})}dt,$ and
that $k=\tau _t$ for $t\in [\sigma _k,\sigma _{k+1}).$ Since $\tau
_T=n,$ we have
\begin{eqnarray*}
\sum_{k=1}^n\frac 1{A_k^{1+ \rho/2}}\, G\left( T_{k-1}\right)
\Delta V_k &=&\sum_{k=1}^n\int_{[\sigma _k,\sigma _{k+1})}\frac 1{A_{\tau
_t}^{1+ \rho/2}}\, G\left( T_{ \tau_t-1}  \right) dt \\
&=&\int_0^T\frac 1{A_{\tau _t}^{1+ \rho/2}}\, G\left( T_{ \tau_t-1} \right) dt.
\end{eqnarray*}
Set  $a_t=c_{*}^2\gamma ^2+T-t.$ Since $\Delta V_{\tau _t}\leq  \epsilon ^2+2 \delta^2$ (cf.\ Lemma \ref{lemma2s}), we see that
\begin{equation}
t\leq V_{\tau _t}\leq V_{\tau _t-1}+\Delta V_{\tau _t}\leq t+ \epsilon
^2+2 \delta^2,\quad t\in [0,T].  \label{BOUND-V}
\end{equation}
Assume that   $c_{*}\geq 4.$ Then we have
\begin{equation}
\frac 12 \, a_t\leq A_{\tau _t}=c_{*}^2\gamma ^2+T-V_{\tau _t}\leq a_t,\quad
t\in [0,T].  \label{BOUND-A}
\end{equation}
Notice that $G(z)$ is symmetric and is non-increasing in $z\geq 0.$ The last bound
implies that
\begin{equation}
J_1\leq 2^{1+\rho/2}  \, \epsilon^{\rho} \int_0^T\frac 1{a_t^{1+ \rho/2}} \, \mathbf{E}  \bigg[ G\bigg( \frac{
u-X_{\tau _t-1}}{a_t^{1/2}}\bigg) \bigg] dt.  \label{W-0}
\end{equation}
Notice also that $G(z)$ is a symmetric integrable function of bounded variation.
By Lemma \ref{LEMMA-APX-2}, it is easy to see that
\begin{equation}
\mathbf{E}  \bigg[ G \left( \frac{u-X_{\tau _t-1}}{{a_t}^{1/2}}\right)  \bigg] \leq
c_{6}\sup_z\Big| \mathbf{P} (X_{\tau _t-1}\leq z)-\Phi (z)\Big| +c_{7} \sqrt{a_t}.
\label{W-1}
\end{equation}
Since $V_{\tau _t-1}=V_{\tau _t}-\Delta V_{\tau _t}$, $ V_{\tau _t}\geq t$ (cf. (\ref{BOUND-V})) and  $\Delta V_{\tau _t}\leq
\epsilon^2 + 2\,\delta^2$, we get
\begin{equation}
V_n-V_{\tau _t-1}\leq V_n - V_{\tau _t}+\Delta V_{\tau _t} \leq 2(\epsilon^2 + \delta^2)+T-t\leq a_t.  \label{CC-2}
\end{equation}
Thus
\begin{eqnarray*}
\mathbf{E}  \left[(X_n-X_{\tau _t-1})^2|\mathcal{F}_{\tau _t-1} \right] &=&\mathbf{E} \bigg[\sum_{k=\tau_t}^n\mathbf{E}  [\xi _k^2 |\mathcal{F}_{k-1}]
\bigg|\mathcal{F}_{\tau_t -1}\bigg] \\
&= &   \mathbf{E} \big[ \left\langle X\right\rangle _n-\left\langle
X\right\rangle_{\tau _t-1}|\mathcal{F}_{\tau_t -1}\big]\\
&\leq&   \mathbf{E}   [V_n-V_{\tau _t -1}|\mathcal{F}_{\tau_t -1} ]    \\
&\leq&   a_t.
\end{eqnarray*}
Then, by Lemma \ref{LEMMA-APX-1}, we deduce that  for any $t\in [0,T],$%
\begin{equation}
\sup_z\Big| \mathbf{P} (X_{\tau _t-1}\leq z)-\Phi (z)\Big| \leq c_8 \, \sup_z\Big|
\mathbf{P} (X_n\leq z)-\Phi (z)\Big| +c_{9}\sqrt{a_t}.  \label{W-2}
\end{equation}
Combining (\ref{fgfgfdgk}), (\ref{W-0}), (\ref{W-1}) and (\ref{W-2}) together,   we obtain
\begin{equation}
\Big|I_1 \Big| \leq c_{10} \, \epsilon^{ \rho}\int_0^T\frac{dt}{a_t^{1+ \rho/2}}\sup_z\Big| \mathbf{P} (X_n\leq
z)-\Phi (z)\Big| +c_{11}\, \epsilon^\rho \int_0^T\frac{dt}{a_t^{ (1+\rho)/2}}.  \label{COMP-1}
\end{equation}
By some elementary computations, it follows that
\begin{equation}
\int_0^T\frac{dt}{a_t^{1+ \rho/2}} \leq \int_0^T\frac{dt}{(c_{*}^2\gamma ^2+T-t)^{1+ \rho/2}} \leq
\frac {1}{c_{*}^{\rho } \gamma^\rho }    \label{W-3}
\end{equation}
and
\begin{displaymath}
\int_0^T\frac{dt}{a_t^{ (1+\rho)/2}} \leq \left\{ \begin{array}{ll}
c_\rho, & \textrm{\ \ \ if $\rho \in (0, 1),$}\\
c \left| \log
  \epsilon  \right|, & \textrm{\ \ \ if $\rho=1.$}
\end{array} \right.
\end{displaymath}
Thus
\begin{equation}
\left| I_1\right| \leq  \frac{c_{12} }{c_{*}^{\rho }}\sup_z\Big| \mathbf{P}(X_n\leq
z)-\Phi (z)\Big| +c_{\rho,1 }\, \widehat{\epsilon},
 \label{W-4}
\end{equation}
where
\begin{displaymath}
\widehat{\epsilon}  = \left\{ \begin{array}{ll}
\epsilon^\rho  +\delta, & \textrm{\ \ \ if $\rho \in (0, 1),$}\\
\epsilon\left| \log   \epsilon  \right|+ \delta , & \textrm{\ \ \ if $\rho=1.$}
\end{array} \right.
\end{displaymath}

\emph{\textbf{b)} Control of} $I_2.$
Note that $0\leq \Delta V_k-\Delta \left\langle X\right\rangle _k\leq
2\delta ^2\mathbf{1}_{\{k=n\}}$. We have
\begin{eqnarray*}
\left|
I_2\right| \leq \mathbf{E}   \Big[ \frac 1{2A_n}\left| \varphi ^{\prime }\left(
T_{n-1}\right) \left( \Delta V_n-\Delta
\left\langle X\right\rangle _n\right) \right|  \Big].
\end{eqnarray*}
Set $\widetilde{G}(z)=\sup_{\left| z-t\right| \leq 1}  |\varphi'(t)|.$ Then $\left| \varphi ^{\prime }(z)\right| \leq \widetilde{G}(z)$ for any real $z.$
Note that  $A_n = c_{*}^2\gamma ^2.$   Then  we get the following estimation:
\[
\left|
I_2\right|\leq \frac{1}{c_{*}^2} \ \mathbf{E}  [ \widetilde{G}\left( T_{n-1}
\right) ].
\]
Notice that $\widetilde{G}(z)$ is   non-increasing in $z\geq 0$, and thus it has  bounded variation on $\mathbf{\mathbf{R}}.$
By Lemmas \ref{lemma2s} and \ref{LEMMA-APX-2}, we obtain
\begin{equation}
\left| I_2\right| \leq \frac{c_{13}}{c_{*}^{2}}\sup_z\Big| \mathbf{P} (X_n\leq z)-\Phi
(z)\Big| +c_{\rho, 2}\,\widehat{\epsilon} .  \label{F-2}
\end{equation}

\emph{\textbf{c)} Control of} $I_3.$  By a two-term Taylor expansion, it follows that
\[
I_3=\frac{1}{8}\,\mathbf{E}   \bigg[\sum_{k=1}^n\frac 1{(A_k-\vartheta _k\Delta A_k)^2}\varphi
^{\prime \prime \prime }\left( \frac{u-X_{k-1}}{\sqrt{A_k-\vartheta _k\Delta A_k} }
\right) \Delta A_k^2\bigg].
\]
Since   $c_{*}\geq 4,\ \Delta A_k \leq 0 $ and
$\left| \Delta A_k\right| =  \Delta V_k  \leq \epsilon^2 + 2\,\delta^2$, we have
\begin{equation}
A_k\leq A_k-\vartheta _k\Delta A_k\leq c_{*}^2\gamma ^2+T-V_k+\epsilon^2 + 2\,\delta^2\leq
2A_k. \label{eqsfion}
\end{equation}
Set $\widehat{G}(z)= \sup_{|t-z|\leq 2}\left| \varphi ^{\prime \prime \prime
}(t)\right|. $ Then
$\widehat{G}(z)$ is symmetric, and is non-increasing in $z\geq 0.$
By (\ref{eqsfion}),
 we obtain
\[
\left| I_3\right| \leq  (\epsilon^2 + 2\,\delta^2)\mathbf{E}  \bigg[ \sum_{k=1}^n\frac 1{A_k^2} \ \widehat{G}\left(
\frac{T_{k-1}}{\sqrt{2 }}\right) \Delta V_k \bigg].
\]
By an argument similar to the proof of (\ref{W-4}),   we get
\begin{eqnarray}
\left| I_3\right| &\leq& \frac{\epsilon^2 + 2\,\delta^2}{c_{*} \gamma }\sup_z\Big| \mathbf{P} (X_n\leq z)-\Phi
(z)\Big| +c_{\rho, 3}\, \widehat{\epsilon}  \nonumber \\
&\leq& \frac{2}{c_{*}   }\sup_z\Big| \mathbf{P} (X_n\leq z)-\Phi
(z)\Big| +c_{\rho, 3}\, \widehat{\epsilon}.  \label{F-3}
\end{eqnarray}

 From  (\ref{lla}),  using (\ref{W-4}), (\ref
{F-2}) and (\ref{F-3}), we have
\begin{eqnarray*}
 \Big| \mathbf{E}  [\Phi _u(X_n,A_n)]-\mathbf{E}  [\Phi _u(X_0,A_0)] \Big|   \leq \frac{c_{14}}{c_{*}^{\rho }}\sup_z\Big| \mathbf{P} (X_n\leq z)-\Phi
(z)\Big| +c_{\rho,4}\, \widehat{\epsilon}.
\end{eqnarray*}
Implementing the last bound in (\ref{llaa022}), we deduce that
\[
\sup_z\Big| \mathbf{P} (X_n\leq z)-\Phi (z)\Big| \leq \frac{c_{15}}{c_{*}^{\rho }}\sup_z\Big|
\mathbf{P} (X_n\leq z)-\Phi (z)\Big| +c_{\rho,5}\,\widehat{\epsilon},
\]
from which, choosing $c_{*}^{\rho }=\max\{ 2c_{15}, 4^\rho\}$, we get
\begin{equation}
\sup_z\Big| \mathbf{P} (X_n\leq z)-\Phi (z)\Big| \leq 2c_{\rho,5} \, \widehat{\epsilon},
\label{F-4}
\end{equation}
which completes the proof of theorem.    \hfill\qed

\subsection{Proof of  Corollary \ref{co20} }
Define $\eta_i=\xi_i$ if $ i \leq v(n),$ $\eta_i=0$ if $ i > v(n).$  Then
 $(\eta_i,\mathcal{F}_i)_{i=0,...,n}$ is also a sequence of martingale differences.  It is easy to see
 that
 \begin{eqnarray*}
\mathbf{E} [| \eta_{i}|^{2+\rho} | \mathcal{F}_{i-1} ] \leq \epsilon^{ \rho} \mathbf{E} [ \eta_{i} ^{2} | \mathcal{F}_{i-1} ].
\end{eqnarray*}
If $v(n)= \sup\{k:\ \langle X\rangle_k \leq1\},$   then
 \begin{eqnarray*}
1 -\mathbf{E} [ \xi_{v(n)+1} ^{2} | \mathcal{F}_{v(n) } ] \leq \sum_{i=1}^n \mathbf{E} [ \eta_{i} ^{2} | \mathcal{F}_{i-1} ]   = \sum_{i=1}^{v(n)} \mathbf{E} [ \xi_{i} ^{2} | \mathcal{F}_{i-1} ]  \leq 1.
\end{eqnarray*}
If $v(n)= \inf\{k:\ \langle X\rangle_k \geq 1\},$ then
 \begin{eqnarray*}
1 \leq \sum_{i=1}^n \mathbf{E} [ \eta_{i} ^{2} | \mathcal{F}_{i-1} ]   = \sum_{i=1}^{v(n)} \mathbf{E} [ \xi_{i} ^{2} | \mathcal{F}_{i-1} ]  \leq 1+\mathbf{E} [ \xi_{v(n)} ^{2} | \mathcal{F}_{v(n)-1} ].
\end{eqnarray*}
Since $\mathbf{E}[\xi_i^2 | \mathcal{F}_{i-1}] \leq \epsilon^2$ for all $i$ (cf.\ Lemma \ref{lemma2s}), we always  have
 \begin{eqnarray*}
\Big|\sum_{i=1}^n \mathbf{E} [ \eta_{i} ^{2} | \mathcal{F}_{i-1} ]-1 \Big| \leq \epsilon ^2.
\end{eqnarray*}
Notice that $\sum_{i=1}^n\eta_{i}= X_{v(n)}.$
Applying Theorem \ref{th0} to $(\eta_i,\mathcal{F}_i)_{i=0,...,n},$ we obtain the desired inequalities.
This completes the proof of  Corollary \ref{co20}.  \qed

\subsection{Proof of Theorem  \ref{th1}.}
To prove Theorem  \ref{th1}, we  use   the following technical  lemma of El Machkouri and Ouchti   \cite{EO07}; see Lemma 1 therein.
\begin{lemma}\label{lemads1}
Let $X$ and $Y$ be  random variables. Then for  any $p\geq 1,$
\begin{eqnarray}
  D(X+Y) \leq 2  D(X) + 3 \Big|\Big|   \mathbf{E} [|Y|^{2p} | X ] \Big|\Big|_1^{1/(2p+1) }.
\end{eqnarray}
\end{lemma}

Following Bolthausen \cite{B82}, consider the stopping time
$$\tau=\sup\{0\leq k\leq n:\ \left\langle X\right\rangle_k \leq 1\}.$$
Assume that $0<\varepsilon  \leq \epsilon. $
Let $r= \lfloor (1- \left\langle X\right\rangle_\tau)/\varepsilon ^2 \rfloor$, where $\lfloor x\rfloor$ stands for the largest integer less than $x$.
Then $r\leq \lfloor 1/\varepsilon ^2 \rfloor $. Let  $N=n+\lfloor 1/\varepsilon ^2 \rfloor +1.$
Let $(\eta_i)_{i\geq 1}$ be a sequence of independent Rademacher random variables, which is also independent of the martingale differences $(\xi_i)_{1\leq i \leq n}$.
For any $i=1,\dots,N,$ define $\xi'_i =\xi_i$ if $i \leq \tau$, $\xi'_i =\varepsilon \eta_i$ if $\tau < i \leq \tau+r,$
$\xi'_{i}=(1-\left\langle X\right\rangle_\tau -r \varepsilon^2 )^{1/2}\eta_i$  if $i=\tau +r +1$, and $\xi'_i=0$ if $\tau +r +1 < i \leq N.$
Clearly, $X'_k=\sum_{i=1}^k \xi'_{i}$, $k=0,\dots,N$ (with $X'_0=0$) is also a martingale sequence  with respect to
the enlarged probability space and the enlarged filtration.
Moreover, it holds  $\left\langle X'\right\rangle_{N}=1$ a.s.\,and condition (A1) is satisfied for $(\xi'_k)_{k=1,\dots,N}.$
Denote by
\begin{displaymath}
\gamma  = \left\{ \begin{array}{ll}
\epsilon^\rho , & \textrm{\ \ \ if $\rho \in (0, 1),$}\\
\epsilon\left| \log   \epsilon  \right|  , & \textrm{\ \ \ if $\rho\geq1.$}
\end{array} \right.
\end{displaymath}
By Theorem  \ref{th1}, it holds, for all $x \in \mathbf{R}$,
\begin{eqnarray}
  \left|\frac{}{} \mathbf{P}(X'_N  \leq x )-  \Phi\left( x\right) \right|   \leq   c_\rho\, \gamma.
\end{eqnarray}
Using Lemma \ref{lemads1}, we get
\begin{eqnarray}
  D(X_n ) &\leq& 2  D(X'_N ) + 3 \Big|\Big| \mathbf{E}[|X_n-X'_N|^{2p} | X'_N ] \Big|\Big|_1^{1/(2p+1) } \nonumber \\
    &\leq& 2\, c_\rho\, \gamma  +  3 \Big( \mathbf{E}[|X_n-X'_N|^{2p} ] \Big)^{1/(2p+1) }. \label{ineq35}
\end{eqnarray}
As $\tau$ is a stoping time, conditionally on $\tau$, the $(\xi_i-\xi_i')_{i\geq\tau+1}$ still forms a martingale difference sequence.
 Using Burkholder's inequality (cf.\ Theorem 2.11 of Hall  and Heyde \cite{HH80}),  we have
 \begin{eqnarray}\label{sfjlmg}
  \mathbf{E}[|X'_N -X_n |^{2p}]  \leq c_p \, \Big(\mathbf{E}\Big[ \Big| \sum_{i=\tau+1}^N\mathbf{E}[(\xi_i-\xi_i')^2 | \mathcal{F}_{i-1} ]\Big|^p \Big] + \mathbf{E}[ \max_{\tau+1 \leq i \leq N} |\xi_i-\xi_i'|^{2p}]\Big) .
 \end{eqnarray}
It is easy to see that
 \begin{eqnarray*}
 \sum_{i=\tau+1}^N\mathbf{E}[(\xi_i-\xi_i')^2 | \mathcal{F}_{i-1} ] &=&  \sum_{i=\tau+1}^n\mathbf{E}[\xi_i^2 | \mathcal{F}_{i-1} ] +\sum_{i=\tau+1}^N\mathbf{E}[\xi'_i\,^2 | \mathcal{F}_{i-1} ] \\
 &=& \langle X\rangle_n+1 - 2\langle X \rangle_\tau .
 \end{eqnarray*}
Notice that $$  1-\mathbf{E}[\xi_{\tau+1}^2 | \mathcal{F}_{\tau}] \leq \langle X\rangle_\tau \leq 1.$$
Hence
 \begin{eqnarray}
 \sum_{i=\tau+1}^N\mathbf{E}[(\xi_i-\xi_i')^2 | \mathcal{F}_{i-1} ] \leq \langle X \rangle_n-1 + 2\mathbf{E}[\xi_{\tau+1}^2 | \mathcal{F}_{\tau}]. \label{ines36}
 \end{eqnarray}
Using the inequality $|a+b|^k \leq 2^{k-1}(|a|^k+ |b|^k), k\geq 1,$ we get
 \begin{eqnarray}
  \mathbf{E}[ \max_{\tau+1 \leq i \leq N} |\xi_i-\xi_i'|^{2p}]  &\leq& 2^{2p-1} \Big( \mathbf{E}[ \max_{\tau+1 \leq i \leq n} |\xi_i|^{2p} ]+ \varepsilon^{2p} \Big) \nonumber \\
  &\leq& 2^{2p-1}  \Big( \mathbf{E}[ \max_{1 \leq i \leq n} |\xi_i|^{2p} ]+  \varepsilon^{2p}  \Big) . \label{ines37}
 \end{eqnarray}
Combining  (\ref{sfjlmg}), (\ref{ines36}) and (\ref{ines37}) together, we deduce that
 \begin{eqnarray}
  \mathbf{E}[|X'_N -X_n |^{2p}]  \leq c_p \, \Big(\mathbf{E}\big[ \big| \langle X\rangle_n-1\big|^p \big] +\mathbf{E}[ \max_{1 \leq i \leq n} |\xi_i|^{2p} ]+  \varepsilon^{2p}\Big) .\nonumber
 \end{eqnarray}
Returning to (\ref{ineq35}) and letting $\varepsilon\rightarrow 0$, we obtain
\begin{eqnarray*}
  D(X_n ) \leq c_{p,\rho} \, \bigg( \gamma+\Big(\mathbf{E}\big[ \big| \langle X\rangle_n-1\big|^p \big] +\mathbf{E}[ \max_{1 \leq i \leq n} |\xi_i|^{2p} ] \Big)^{1/(2p+1) } \bigg).
\end{eqnarray*}
This completes the proof of Theorem   \ref{th1}. \hfill\qed

\section*{Appendix\label{Appendix}}

In the proof of Theorem \ref{th0}, we make use of the following two technical lemmas due to Bolthausen  (cf.\ Lemmas 1 and 2 of \cite{B82}).
\begin{lemma}
\label{LEMMA-APX-1}Let $X$ and $Y$ be random variables. Then
\[
\sup_u\Big| \mathbf{P}\left( X\leq u\right) -\Phi \left( u \right) \Big|
\leq c_1\sup_u\Big| \mathbf{P}\left( X+Y\leq u\right) -\Phi \left( u \right)
\Big| +c_2 \Big\| \mathbf{E}\left[ Y^2|X\right] \Big\| _\infty ^{1/2}.
\]
\end{lemma}
\begin{lemma}
\label{LEMMA-APX-2}Let $G(x)$ be an integrable function on $\mathbf{\mathbf{R}}$ of bounded variation $||G||_V$,
$X$ be a random variable and $a,$ $b\neq 0$ are real numbers. Then
\[
\mathbf{E}\bigg[ \, G\left( \frac{X+a}b\right) \bigg] \leq ||G||_V \sup_u\Big| \mathbf{P}\left( X\leq u\right)
-\Phi \left( u\right) \Big| + ||G||_1 \, |b|,
\]
where  $||G||_1$ is the $L_1(\mathbf{R})$ norm of $G(x).$
\end{lemma}

\section*{Acknowledgements}
The author  would like to thank editor and reviewers for their comments.
This work has been partially supported by the National Natural Science Foundation
of China (Grant nos.\,11601375 and 11626250).

\end{document}